\documentclass[11pt,a4paper]{article}
\usepackage{amsmath,amssymb,latexsym,amsthm}  
\usepackage[english]{babel}
\usepackage[latin2]{inputenc}

\begin{document}

\newtheorem{Dfn}{Definition}
\newtheorem{Theo}{Theorem}
\newtheorem{Lemma}[Theo]{Lemma}
\newtheorem{Prop}[Theo]{Proposition}
\newcommand{\Pro}{\noindent{\em Proof. }}
\newcommand{\Rem}{\noindent{\em Remark. }}
\renewcommand{\thefootnote}{}

\title{Bol loops as sections in  semi-simple Lie groups of small dimension}
\author{\'Agota Figula} 
\date{}
\maketitle

\begin{abstract}  Using the relations between the theory of differentiable 
Bol loops and the theory of affine symmetric spaces  we 
classify 
all connected differentiable Bol loops having 
an at most $9$-dimensional semi-simple Lie group  as the group 
topologically generated by their left translations. We show that all 
these Bol loops are isotopic to direct products of Bruck loops of hyperbolic type or to Scheerer extensions of Lie groups   by Bruck loops of hyperbolic type.   
\end{abstract}
 


\noindent
{\footnotesize {{\em Thanks: } This paper was supported by DAAD.}

\section{Introduction}
\label{intro}

\bigskip
\noindent
In \cite{quasigroups} and 
\cite{loops} the authors have thoroughly studied  the relations between smooth Bol loops and homogeneous spaces of Lie groups. Using 
their point of view  we treat the 
connected  differentiable Bol loops $L$ as images of global differentiable 
sections $\sigma :G/H \to G$, where $G$ is a connected Lie 
group, $H$ is a 
closed subgroup containing no non-trivial normal subgroup of $G$ and  for 
all $r,s \in \sigma (G/H)$ the element $rsr$ lies in 
$\sigma (G/H)$. 
The differentiable Bruck loops, these are Bol loops $L$ satisfying the 
identity $(x y)^{-1}=x^{-1} y^{-1}$, $x,y \in L$, have 
realizations on differentiable affine symmetric spaces $G/H$, where 
$H$ is the set of fixed elements of an involutory automorphism of $G$ and 
$\sigma (G/H)$ is the exponential image of the $(-1)$-eigenspace of the corresponding   automorphism of the Lie algebra ${\bf g}$ of $G$. An 
impotant subclass of Bruck loops are the Bruck loops of hyperbolic 
type which correspond to Lie groups $G$ and involutions $\tau $ fixing 
elementwise a 
maximal compact subgroup of $G$ (cf. \cite{freudenthal}, 64.9, 
64.10). With this notions our main result reads as follows.

\bigskip
\noindent
{\bf Main Theorem} {\it Let $L$ be a connected differentiable Bol loop having an at most $9$-dimensional semi-simple Lie group $G$ as 
the group 
topologically generated by its left translations. Then $L$ is isotopic to a 
direct product of Bruck loops of hyperbolic type or to a Scheerer extension of a Lie group by a Bruck loop of hyperbolic type.  
\newline
If $\hbox{dim}\ G \le 5$ then $L$ is isotopic to the hyperbolic plane loop 
$\mathbb H_2$ and $G$ is isomorphic to $PSL_2(\mathbb R)$ (cf. 
\cite{loops}, Sec. 22).  
\newline
If $\hbox{dim}\ G=6$ then $L$ is isotopic either to the direct product 
$\mathbb H_2 \times \mathbb H_2$ or to the $3$-dimensional hyperbolic 
 space loop $\mathbb H_3$ (cf. \cite{figula1}, p. 446) or to a Scheerer extension 
 of a $3$-dimensional simple Lie group $G_1$ by $\mathbb H_2$. 
 In the first case $G$ is 
isomorphic to $PSL_2(\mathbb R) \times PSL_2(\mathbb R)$, in the second 
case to $PSL_2(\mathbb C)$, in the third case to $PSL_2(\mathbb R) \times G_1$. 
\newline
If $7 \le  \hbox{dim}\ G \le 8$ then either $L$ is isotopic to the complex 
hyperbolic plane loop having $PSU_3(\mathbb C,1)$ as the group 
topologically 
generated by its left translations (cf. \cite{lowen}, Prop. 5.1, p.152) or 
$L$ is isotopic to the $5$-dimensional Bruck loop of 
hyperbolic type having the group $PSL_3(\mathbb R)$ as the group 
topologically generated by its left translations. 
\newline
If $\hbox{dim}\ G=9$ then $L$ is isotopic either to 
$\mathbb H_2 \times \mathbb H_2 \times \mathbb H_2$ or to 
$\mathbb H_2 \times \mathbb H_3$ or to a Scheerer extension of a $6$-dimensional semi-simple Lie group $G_2$ by $\mathbb H_2$ or to a Scheerer extension of a $3$-dimensional simple Lie group $G_1$ by $\mathbb H_2 \times \mathbb H_2$ respectively by $\mathbb H_3$.  
In the first case 
$G$ is isomorphic to 
$PSL_2(\mathbb R) \times PSL_2(\mathbb R) \times PSL_2(\mathbb R)$, 
in the second case to 
$PSL_2(\mathbb R) \times PSL_2(\mathbb C)$, in the third case to $PSL_2(\mathbb R) \times G_2$ and in the fourth case to $PSL_2(\mathbb R) \times PSL_2(\mathbb R) \times G_1$ respectively to 
$PSL_2(\mathbb C) \times G_1$. }

\medskip
\noindent  
All known differentiable connected Bol loops having a semi-simple 
Lie group as the group topologically generated by 
their left translations  are isotopic to  direct products of 
Bruck loops of hyperbolic type 
or to Scheerer extensions 
of Lie groups by Bruck loops of hyperbolic type since they are 
constructed using Cartan involutions.    

\medskip
\noindent
For the classification of differentiable  Bol loops $L$ up to isotopisms 
we  proceed in the following way (\cite{quasigroups}, pp. 424-425, \cite{loops}, pp. 78-79 and Proposition 1.6, p. 18).  In  the Lie algebra ${\bf g}$ of the 
group $G$ topologically generated by the left translations of $L$  we 
determine the $(-1)$-eigenspaces ${\bf m}$ for 
 all involutory automorphisms of ${\bf g}$. 
After this we seek for any ${\bf m}$ a system of representatives ${\bf h}^*$ 
of the sets $\{Ad_g {\bf h}; g \in G \}$ which consists of subalgebras 
with ${\bf g}={\bf m} \oplus Ad_g {\bf h}={\bf m} \oplus {\bf h}^*$. 
Any triple $(G,\exp {\bf h}^*, \exp {\bf m})$ determines a local Bol 
loop. Global differentiable Bol loops $L$ correspond precisely to those 
exponential images of ${\bf m}$, 
which form a system of representatives for the cosets of $\exp {\bf h}$ in 
$G$. To show which local Bol loop is extendible to a global one we need a 
futher analytic treatment since there are much more local than 
global differentiable Bol loops.   

\medskip
\noindent
I would thank to the referee for his useful remarks.

\section{Some basic notions of the theory of Bol loops} 
\label{sec:1}

\bigskip
\noindent 
A binary system  $(L, \cdot )$ is called a loop if there exists an element 
$e \in L$ such that $x=e \cdot x=x \cdot e$ holds for all $x \in L$ and the 
equations 
$a \cdot y=b$ and $x \cdot a=b$ have precisely one solution which we denote 
by $y=a \backslash b$ and $x=b/a$. Two loops $(L_1, \circ )$ and 
$(L_2, \ast )$ are 
called isotopic if there are three bijections 
$\alpha ,\beta ,\gamma : L_1 \to L_2 $ such that 
$\alpha (x) \ast \beta (y)=\gamma (x \circ y)$ holds for any 
$x,y \in L_1$. Isotopism of loops is an equivalence relation. Let 
$(L_1, \cdot )$ and $(L_2, \ast )$ be two loops. The set  
$L=L_1 \times L_2= \{ (a,b) \ |\ a \in L_1, b \in L_2 \}$ with the 
componentwise multiplication  is again a loop, which is called the direct 
product of $L_1$ and 
$L_2$, and the loops $(L_1, \cdot )$,  $(L_2, \ast )$ are subloops of $L$.   

\medskip
\noindent
A loop $L$ is called a Bol loop if for any two left translations 
$\lambda _a, \lambda _b$ the product $\lambda _a \lambda _b \lambda _a $ 
is 
again a left translation of $L$. If $L_1$ and $L_2$ are Bol loops, then the 
direct product $L_1 \times L_2$ is again a Bol loop. Every subloop of a Bol 
loop satisfies the Bol identity.

\medskip
\noindent
The theory of differentiable loops $L$ is essentially the theory of the 
smooth binary operations  
$(x,y) \mapsto x \cdot y, \ (x,y) \mapsto x/y, \ (x,y) \mapsto x 
\backslash y $ on the connected differentiable manifold $L$. If  $L$ is a 
connected differentiable Bol loop then the left translations 
$\lambda _a= y \mapsto a \cdot y :L \to L$, $a \in L$, which are 
diffeomorphisms  of $L$, topologically generate  a connected Lie  group 
$G$ (cf. \cite{loops}, 
p. 33; \cite{quasigroups}, pp. 414-416). Moreover the manifold $L$ is 
parallelizable since the set of the left translations is sharply transitive.

\medskip
\noindent
Every connected differentiable Bol loop is isomorphic to a  loop $L$ 
realized on the homogeneous space $G/H$, 
where $G$ is a connected Lie group, $H$ is a connected 
closed subgroup which is not allowed  to  contain a  non-trivial  normal 
subgroup  of $G$ and 
$\sigma : G/H \to G$ is  a differentiable 
section with $\sigma (H)=1 \in G$ such that the 
subset $\sigma (G/H)$ generates $G$  and for all $r,s \in \sigma (G/H)$ the 
element $rsr$ is contained in 
$\sigma (G/H)$ (cf. \cite{loops}, p. 18 and Lemma 1.3, 
p. 17, \cite{kiechle}, Corollary 3.11, p. 51). The multiplication  of  $L$ 
on the  
space $G/H$  is  defined by $x H \ast y H=\sigma (x H) y H$ and  the group 
$G$ is  the group topologically generated by the left 
translations of $L$.  
\newline 
Let ${\bf m}=T_1 \sigma (G/H)$ be the tangent space of $\sigma (G/H)$ at 
$1 \in G$. If $(\bf{g, [.,.]})$ 
respectively ${\bf h}$ denotes the Lie algebra of $G$ respectively of $H$ 
then one has 
$\bf{g}=\bf{m} \oplus \bf{h}$, 
$\big[ [\bf{m},\bf{m}], \bf{m} \big] \subseteq \bf{m}$ and ${\bf m}$ 
generates ${\bf g}$. 
>From every Bol triple $(({\bf g}, [.,.]),\ {\bf h},\ {\bf m})$ we can 
obtain with a canonical construction  a triple 
$({({\bf g}^*, [.,.]^*), {\bf h}^*, {\bf m}^*})$, 
such that ${\bf m}^*$ is the $(-1)$-eigenspace of an involutory automorphism of ${\bf g}^*$ and the subalgebra ${\bf h}^*$ complements ${\bf m}^*$ in 
${\bf g}^*$ 
(cf. \cite{quasigroups}, 
pp. 424-425). If the Lie algebra ${\bf g}$ is semi-simple  then it is 
isomorphic to ${\bf g^*}$ and it is the Lie algebra of the displacement group for the 
symmetric space belonging to the Lie triple system $A:=\big( {\bf m},(.,.,.) \big)$, with $(.,.,.)=\big[ [.,.], .\big]$ (cf. \cite{loos}, pp. 78-79).  
Hence one has 
${\bf g}={\bf m} \oplus [ {\bf m}, {\bf m}]$ (cf. \cite{loops}, Section 6). 

\medskip
\noindent
\Rem Let ${\bf g}$ be the Lie algebra of the group topologically generated 
by 
the left translations of a connected 
differentiable Bol loop $L$ such that  ${\bf g}$ is  the direct sum  
${\bf g}_1 \oplus {\bf g}_2$ of  
Lie algebras ${\bf g}_i$, $i=1,2$. If an involutory automorphism $\tau $ of 
${\bf g}$ has the shape 
$(\tau _1, id)$, where $\tau_1 $ is an  involutory automorphism of 
${\bf g}_1$ and 
$id: {\bf g}_2 \to {\bf g}_2$ is the identity map, then one has 
${\bf m}={\bf m}_1 \oplus {\bf g}_2$ and 
${\bf h}={\bf h}_1 \oplus \{0 \}$, where 
${\bf m}_1$ is the $(-1)$-eigenspace  and ${\bf h}_1$ the $(+1)$-eigenspace 
of $\tau_1$. The Bol loop 
$L$ corresponding to such a triple $({\bf g}, {\bf h}, {\bf m})$ is the 
direct product of a Bol loop $\tilde{L}$ 
isotopic to the Bruck loop realized on the symmetric space $\exp {\bf g}_1 / \exp {\bf h}_1$ and  the Lie group $\exp {\bf g}_2$. 

\medskip
\noindent
Let $L_1$ be a loop defined on the factor space $G_1 /H_1$ with respect to a section $\sigma _1:G_1 /H_1 \to G_1$ the image of which is the set $M_1 \subset G_1$. Let $G_2$ be a group, let $\varphi :H_1 \to G_2$ be a homomorphism and 
$(H_1, \varphi (H_1))=\{ (x, \varphi (x))\ |\ x \in H_1 \}$.  A loop $L$ is called a 
Scheerer extension of $G_2$ by $L_1$ if the loop $L$ is defined on the factor space 
$(G_1 \times G_2)/(H_1, \varphi (H_1))$ with respect to the section $\sigma :(G_1 \times G_2)/(H_1, \varphi (H_1)) \to G_1 \times G_2$ the image of which is the set $M_1 \times G_2$.

\medskip
\noindent
\begin{Lemma} Let ${\bf g}$ be the Lie algebra of the group topologically 
generated by the left translations of a  connected differentiable Bol loop 
$L$ such that $\hbox{dim}\ L \ge 4$ and the tangent space $T_1 L$ does not 
contain a Lie algebra as a direct factor. 
Let  ${\bf g}$ be the direct product of  simple Lie algebras with 
 $\hbox{dim}\ {\bf g} \le 9$.  
\newline
i) If $\hbox{dim} \ {\bf g}=6$ then $\hbox{dim} \ L=4$.
\newline
ii) If $\hbox{dim} \ {\bf g}=9$ then $\hbox{dim} \ L \in \{ 5,6 \}$. 
\end{Lemma}

\noindent
\Pro Let $\tau $ be the involutory automorphism corresponding to $L$.  
Neither the $(-1)$-eigenspace nor the $(+1)$-eigenspace of $\tau $ can 
contain a simple direct factor of  
${\bf g}$. The $(-1)$-eigenspace of each involutory automorphism of a 
$3$-dimensional simple Lie algebra is 
$2$-dimensional (cf. \cite{furness}, pp. 44-45). Hence if 
${\bf g}={\bf g}_1 \oplus {\bf g}_2$ with simple 
$3$-dimensional Lie algebras ${\bf g}_i$ then $\hbox{dim}\ L =4$. If  
${\bf g}= {\bf g}_1 \oplus {\bf g}_2 
\oplus {\bf g}_3$ with $3$-dimensional simple Lie algebras ${\bf g}_i$, 
$(i=1,2,3)$ then the $(+1)$-eigenspace of an involutory automorphism of 
${\bf g}$ has  dimension either $3$ or $4$.  
 Let 
${\bf g}={\bf g}_1 \oplus {\bf g}_2$, where ${\bf g}_i$ are simple 
Lie algebras with $\hbox{dim} \ {\bf g}_1=6$ and $\hbox{dim} \ {\bf g}_2=3$. The dimension of a 
$(-1)$-eigenspace of an involutory automorphism of ${\bf g}_1$ is either 
$3$ or $4$ (cf. \cite{lowen}, p. 153). 
Hence $\hbox{dim} \ L \notin \{4,7,8 \}$.   
\qed

\medskip
\noindent
>From topological reasons we obtain 

\begin{Lemma}  Let $G$ be isomorphic to the Lie group  
$G_1 \times G_2$, such that 
$G_2 \cong SO_3(\mathbb R)$ and for the subgroup $H$ of $G$ one has  
$H=H_1 \times H_2$ with $1 \neq H_2 \le G_2$. 
Then $G$ cannot be the group topologically generated by the left 
translations of a topological loop. 
\end{Lemma}

\noindent
\Pro The  factor space $G/H$ is a  
 topological product of spaces having as a factor the $2$-sphere or the 
projective plane, which are not parallelizable. 
\qed

\begin{Lemma} Let $G$ be a Lie group isomorphic to $K \times SO_3(\mathbb R) \times SO_3(\mathbb R)$ and let $H$ 
be a subgroup of $G$ such that 
$H=H_1 \times \{ (a,a) \ | \ a \in SO_3(\mathbb R) \}$  with $H_1 \le K$. 
Then there is no Bol loop $L$ corresponding to the pair $(G,H)$. 
\end{Lemma}

\noindent
\Pro The loop $L$ would be a product of a loop $L_1$ corresponding to the 
factor space $K/H_1$ with a compact 
proper loop $L_2$ having $SO_3(\mathbb R) \times SO_3(\mathbb R)$ as the 
group topologically generated by 
its left translations (cf. \cite{loops}, Proposition 1.18, p. 26). But such 
a loop $L_2$ does not exist 
(cf. \cite{loops}, Corollary 16.9, p. 204). 
\qed

\begin{Prop} There is no connected differentiable  Bol loop such that the 
group $G$ topologically generated by 
its left translations is a  compact Lie group $G$ with 
$\hbox{dim} \ G \le 9$. 
\end{Prop} 
\Pro If $G$ is a  quasi-simple Lie group then it admits a continuous section if and only if $G$ is 
locally isomorphic to  $SO_8(\mathbb R)$ (cf. \cite{scheerer}, pp. 149-150).  Since 
$\hbox{dim} \ L \ge 4$ (\cite{loops}, 
Corollary 16.8, p. 204) it follows from \cite{loops}, Theorem 16.1 that $G$ 
is locally isomorphic to 
$SO_3(\mathbb R) \times SO_3(\mathbb R) \times SO_3(\mathbb R)$. This is 
excluded by Theorem 16.7 in \cite{loops}, p. 198 and Lemmata 2 and 3.  
\qed 

\medskip
\noindent
An important tool to eliminate  certain stabilizers $H$ is the fundamental 
group $\pi _1$ of a connected topological 
space.  

\begin{Lemma} Denote by $G$  a connected Lie group and by 
 $H$  a connected  
subgroup of $G$. Let 
$\sigma : G/H \to G$ be a global section. Then 
$\pi_1(K) \cong \pi_1(\sigma (G/H)) \times \pi_1(K_1)$, where 
$K$ respectively $K_1$ is  a maximal compact subgroup of $G$ respectively of 
$H$. 
\end{Lemma}

\noindent
\Pro Since $G$ respectively $H$ is homeomorphic to a topological product 
$K \times \mathbb R^n$ respectively to 
$K_1 \times \mathbb R^n$ (cf. \cite{montgomery}, p. 178) one has 
$\pi_1(G)=\pi_1(K)$ and $\pi_1(H)=\pi_1(K_1)$. 
The group $G$ is homeomorphic to the topological product 
$\sigma (G/H) \times H$. Hence 
$\pi_1(G)=\pi_1(K) \cong 
\pi_1(\sigma(G/H)) \times H) \cong \pi_1(\sigma (G/H)) \times \pi_1(K_1)$ 
(cf. \cite{hu}, Theorem 2.1, p. 144).  
\qed 

\medskip
\noindent
Since the Bol loops are strongly left 
alternative (Definition 5.3. in \cite{loops}) every global Bol loop $L$ 
contains the 
exponential image of the tangent  space  ${\bf m}$ at  $e \in L$. 
In the  discussion 
which  submanifolds  $\exp {\bf m}$ can be extended to a global section we 
use the following Lemma of  \cite{figula1}.

\begin{Lemma} Let $L$ be a differentiable loop and denote by ${\bf m}$ the 
tangent space $T_1 \sigma(G/H)$, where 
$\sigma : G/H \to G$ is the section corresponding 
to $L$. Then ${\bf m}$ does not contain any element of $Ad_g {\bf h}$ for 
some $g \in G$.  Moreover, every element 
of $G$ can be uniquely written as a product of an 
element of $\sigma(G/H)$ with an element of $H$. 
\end{Lemma}

\begin{Lemma} Let $\sigma :G/H \to G$ be a continuous section, where  
$G \cong PSL_2(\mathbb R)$ and $H$ is the group 
 $\left \{ \left ( \begin{array}{ll}
l & b \\
0 & l^{-1} \end{array} \right );\ b \in \mathbb R \right \}$  such that 
either $l=1$ for all $b \in \mathbb R$ or 
$0 < l \in \mathbb R$. Then the image $\sigma (G/H)$ cannot contain the 
manifold $M=\left \{ \left ( \begin{array}{cc} 
x+y & z \\
z & x-y \end{array} \right ); \  y,z \in \mathbb R, x \ge 1, x^2-y^2-z^2=1 
\right \}$.   
\end{Lemma}

\noindent
\Pro The coset 
$g(c) H=\left( \begin{array}{cc} 
1+c & 1 \\ 
c   & 1 \end{array} \right) H,
$  
$c > -1$,
contains the element 
$
s(c)= \left( \begin{array}{cc} 
1+c & c \\ 
c   &  \frac{c^2+1}{1+c} \end{array} \right) \in M. 
$ 
But 
$\lim \limits_{c \to -1} \sigma (g(c) H)$ should be  
$\lim \limits_{c \to -1} s(c)$ which is a contradiction. 
\qed

\medskip
\noindent
\Rem If the Lie group $G$ is the group topologically generated by the left 
translations of two Bol loops $L_1$ and $L_2$ such that the corresponding 
symmetric spaces are isomorphic then $L_1$ and $L_2$ are isotopic (cf. \cite{loops}, Theorem 1.11, pp. 21-22).

\medskip
\noindent
In our computations  we often use the following facts. 

\medskip
\noindent
As a real basis of $\mathfrak{sl_2}(\mathbb R)$ we choose  

\medskip
\noindent
(1) \centerline{$e_1=\left (\begin{array}{rr} 
1 & 0 \\ 
0 & -1 \end{array} \right )$, \ $e_2=\left (\begin{array}{cc} 
0 & 1 \\
1 & 0 \end{array} \right )$, \ $e_3=\left (\begin{array}{rr}
0 & 1 \\
-1 & 0 \end{array} \right )$} 
(cf. \cite{hilgert}, pp. 19-20).  

\medskip
\noindent
Then the Lie algebra multiplication is given by 

\medskip
\noindent
\centerline{$[e_1,e_2]=2 e_3$, \ $[e_1,e_3]=2 e_2$, \ $[e_3,e_2]=2 e_1$. }

\medskip
\noindent
An element $X=\lambda _1 e_1+ \lambda _2 e_2+ \lambda _3 e_3 
\in \mathfrak{sl_2}(\mathbb R)$ is elliptic, parabolic or 
hyperbolic according to whether 

\medskip
\noindent
\centerline{$k(X,X)=\lambda _1^2+\lambda _2^2-\lambda _3^2$\ \  is smaller, 
equal, or greater $0$. }

\medskip
\noindent
{\bf 1.1} The basis elements $e_1,e_2$ are hyperbolic, $e_3$ 
is elliptic and the elements $e_2+e_3$, $e_1+e_3$ are both parabolic. 
All elliptic elements are conjugate under $Ad\ G$ to $e_3$, 
 all hyperbolic elements to $e_1$ and 
all parabolic 
elements  to $e_2+e_3$ (cf. \cite{hilgert}, p. 23). 
There is precisely one conjugacy class $\mathcal{C}$ of the two dimensional 
subgroups of $PSL_2(\mathbb R)$; as a representative of $\mathcal{C}$ we 
choose 

\medskip
\noindent
\centerline{${\mathcal L}_2= \left \{ \left ( \begin{array}{ll}
a & b \\
0 & a^{-1} \end{array} \right );\ a>0, b \in \mathbb R \right \}$. }
 
\medskip
\noindent
The Lie algebra of ${\mathcal L}_2$ is generated by the elements $e_1$, 
$e_2+e_3$. 

\medskip
\noindent
{\bf 1.2} As a real basis of the Lie algebra 
$\mathfrak{so_3}(\mathbb R) \cong \mathfrak{su_2}(\mathbb C)$ we may  choose 
the basis elements $\{\hbox{i} e_1, \hbox{i} e_2, e_3\}$, where 
$\hbox{i}^2=-1$. 
Every elements of $\mathfrak{so_3}(\mathbb R)$ is conjugate to $e_3$.

\section{ Bruck loops of hyperbolic type}  
\label{sec:2}

\bigskip
\noindent
Now we give a procedure to construct many Bruck loops having a non-compact 
Lie group as the group topologically generated by the left translations. 

\begin{Theo}
Let $G$ be a simple non-compact Lie group, let $H$ be a  maximal compact 
subgroup of $G$ and let $\tau $ be the involutory automorphism of the Lie 
algebra ${\bf g}$ of $G$ such that  the Lie algebra ${\bf h}$ of $H$ is the 
$(+1)$-eigenspace of $\tau $.  

\medskip
\noindent 
a) The factor space $G/H$ is a Riemannian symmetric space diffeomorphic to 
 the manifold $\exp {\bf m}$, where ${\bf m}$ is the $(-1)$-eigenspace of 
$\tau $.  The group $G$ is the group 
of displacements of $G/H$ and $\exp {\bf m}=\{ \sigma _{xH} \sigma _H \}$, 
where $ \sigma _{xH}$ is the reflection at the point $xH$. Any coset $xH$ 
contains precisely one element of $\exp {\bf m}$.

\medskip
\noindent 
b) The section $\sigma :G/H \to G$ assigning to the coset $xH$ the element 
of $\exp {\bf m}$ contained in $xH$ defines a global Bruck loop $L$ on $G/H$ 
by $(xH) \ast (yH)=\sigma (G/H) yH$. 
\end{Theo}

\noindent
\Pro According to \cite{loos}, Proposition 1.7, p. 148, the factor space 
$G/H$ is a Riemannian symmetric  space having the group $G$ as its group of 
displacements, such that $\exp {\bf m}$ consists of the products of the 
reflection at $H$ and a reflection at arbitrary point  
(cf. \cite{loos}, p. 64 and Proposition 3.2, p. 95 and Theorem 1.3, p. 73). 
 It follows 
from \cite{loos}, Theorem 3.2, p. 165, that $G/H$ is diffeomorphic to 
$\exp {\bf m}$,  any coset $xH$ contains 
precisely one element of $\exp {\bf m}$ and  the section 
$\sigma :G/H \to G$ is differentiable.  For elements 
 $r=\sigma _{xH} \sigma _H$, 
$s=\sigma _{yH} \sigma _H$ of  $\exp {\bf m}$ one has 
$rsr=( \sigma _{xH} \sigma _H) ( \sigma _{yH} \sigma _H) 
( \sigma _{xH} \sigma _H)= (\sigma _{xH}( \sigma _H \sigma _{yH} \sigma _H) 
\sigma _{xH}) \sigma _H \in \exp {\bf m}$. Hence the multiplication $\ast :
G/H \times  G/H \to G/H$ defines on $G/H$ a global differentiable Bruck loop 
$L$ (cf. \cite{loops}, Proposition 9.25, p. 118). 
\qed 

\medskip
\noindent
Another proof of this theorem is given in 
\cite{krammer}, Theorem 3.3, p. 319.

\bigskip
\noindent
Any simple non-compact Lie group $G$  admits an involutory automorphism 
$\tau $ 
the centralizer of which is a maximal compact subgroup of $G$. This 
Cartan involution  
$\tau $ determines a  symmetric space ${\cal S}(\tau )$ of hyperbolic type 
(cf. \cite{freudenthal}, 64.9, p. 375). For this reason we call the Bruck 
loop 
realized on ${\cal S}(\tau )$ a differentiable Bruck loop of hyperbolic 
type.

\section{ Bol loops corresponding to  simple Lie groups}
\label{sec:3}

\bigskip
\noindent
First we investigate the Lie groups locally isomorphic to 
$PSL_2(\mathbb C)$. A 
real basis of the Lie algebra $\mathfrak{sl_2}(\mathbb C)$ is given by  
$\{ e_1, e_2, e_3, \hbox{i} e_1, \hbox{i} e_2, \hbox{i} e_3 \}$, where
$\{e_1, e_2, e_3 \}$ is the basis of $\mathfrak{sl_2}(\mathbb R)$ described 
by (1). 

\medskip
\noindent
According to \cite{lowen}, p. 153, there is only one conjugacy class of 
involutory 
automorphisms of $G=PSL_2(\mathbb C)$ leaving an at most $2$-dimensional 
subgroup 
$H$ of $G$ elementwise fixed.  A representative of this class is the map 
$\tau :PSL_2(\mathbb C) \to PSL_2(\mathbb C): X \mapsto A^{-1} X A$  with  
$A=\pm 
\left ( \begin{array}{rr} 
0 & 1 \\ 
-1 & 0 \end{array} \right )$. The centralizer of $\tau $ is the group 
$H=\left \{ \pm \left (\begin{array}{rr} 
a & b \\
-b & a \end{array} \right ),\ a, b \in \mathbb C, a^2+b^2=1 \right \}$.   
\newline
The 
Lie algebra ${\bf h}$ of  $H$   is generated by 
$ e_3,  \hbox{i} e_3 $ and the tangent space  ${\bf m}$ of the 
corresponding $4$-dimensional symmetric space has as basis elements 
$e_1, e_2,  \hbox{i} e_1,  \hbox{i} e_2$. 

\begin{Lemma} Let ${\bf g}$ be the Lie algebra $\mathfrak{sl_2}(\mathbb C)$.
Any $2$-dimensional subalgebra ${\bf h}$ of ${\bf g}$  
has (up to conjugation) one 
of the following forms: 

\medskip
\noindent 
 \ \ \  ${\bf h}_1=\langle e_3, \hbox{i} e_3 \rangle $, \ \  \  
${\bf h}_2=\langle e_1, e_2+e_3 \rangle $, \ \  \  
${\bf h}_3=\langle \hbox{i} (e_2+e_3), e_2+e_3 \rangle $. 

\medskip
\noindent 
The element  $e_2+e_3 \in {\bf h}_2$ as well as 
 $ \hbox{i} (e_2+e_3) \in {\bf h}_3$  is conjugate under $Ad \ G$  to 
$e_1+\hbox{i} e_2 \in {\bf m}$. 
The element  $e_3 \in {\bf h}_1$ respectively $\hbox{i} e_3 \in {\bf h}_1$ is 
conjugate 
to $\hbox{i} e_1 \in {\bf m}$ respectively to  $-e_1 \in {\bf m}$.
\end{Lemma}

\noindent
\Pro The first assertion follows from Theorem 15 in \cite{lie}, p. 129. 
\newline
Futhermore one has 
$Ad_{g_1}(e_2+e_3)=e_1+\hbox{i} e_2$, 
 $Ad_{g_2}(\hbox{i}(e_2+e_3))=e_1+\hbox{i} e_2$,  
$Ad_{g_3}(e_3)=\hbox{i} e_1$,  
$Ad_{g_3}(\hbox{i} e_3)=- e_1$, where  
$g_1=\pm \left (  \begin{array}{cc}
-1+\hbox{i} & 0 \\
-\frac{1}{2}+\frac{1}{2} \hbox{i} & -\frac{1}{2}-\frac{1}{2} \hbox{i} 
\end{array} \right ) $, 
$g_2=\pm \left (  \begin{array}{cc}
- \frac{2}{\sqrt 2} & 0 \\
\frac{1}{\sqrt 2} \hbox{i} & -\frac{1}{\sqrt 2} \end{array} \right ) $ and 
 $g_3=\pm \left (  \begin{array}{cc}
1  & \frac{1}{2} \hbox{i} \\
\hbox{i} & \frac{1}{2} \end{array} \right )$. 
\qed

\medskip
\noindent
Lemma 9  yields the    

\begin{Prop} There is no at least $4$-dimensional Bol loop having a group 
locally 
isomorphic to $PSL_2(\mathbb C)$ 
as the group topologically generated by its left translations. 
\end{Prop}

\bigskip 
\noindent
Now let $G$ be locally isomorphic to the non-compact Lie group 
$PSU_3(\mathbb C,1)$. 
The Lie algebra 
${\bf g}=\mathfrak{su_3}(\mathbb C,1)$ of $G$ can be treated as the Lie 
algebra of matrices 
\[(\lambda _1 e_1+\lambda _2 e_2+\lambda _3 e_3+\lambda _4 e_4+\lambda _5 e_5+\lambda _6 e_6+\lambda _7 e_7+\lambda _8 e_8) \mapsto \]
\[ \left ( \begin{array}{ccc}
- \lambda _1 \hbox{i}  & - \lambda _2 -\lambda _3 \hbox{i} & \lambda_4+ \lambda _5 \hbox{i} \\
 \lambda _2 -\lambda _3 \hbox{i} & \lambda _1 \hbox{i}+ \lambda _6 \hbox{i} & \lambda _7 +\lambda _8 \hbox{i} \\
 \lambda _4 -\lambda _5 \hbox{i} & \lambda _7 - \lambda_8 \hbox{i}  & - \lambda _6 \hbox{i} \end{array} \right ); \lambda_j \in 
\mathbb R, j=1,\cdots ,8. \] 
Then the  multiplication  in ${\bf g}$ is given by the following: 
\[ [e_1,e_6]=0,\  [e_3,e_2]= 2 e_1,\  [e_4,e_5]= 2(e_1-e_6),\  [e_8,e_7]=2 e_6, \]
\[ [e_6,e_3]=[e_7,e_4]=[e_8,e_5]=\frac{1}{2} [e_1,e_3]= e_2, \]
\[ [e_2,e_6]=[e_4,e_8]=[e_7,e_5]=\frac{1}{2} [e_2,e_1]=e_3, \]
\[ [e_7,e_2]=[e_3,e_8]=[e_5,e_6]= [e_1,e_5]=e_4, \] 
\[ [e_8,e_2]=[e_7,e_3]=[e_6,e_4]= [e_4,e_1]=e_5, \]
\[ [e_2,e_4]=[e_3,e_5]=[e_8,e_1]=\frac{1}{2} [e_8,e_6]=e_7, \] 
\[ [e_2,e_5]=[e_4,e_3]=[e_1,e_7]=\frac{1}{2} [e_6,e_7]=e_8.  \]
There are two conjugacy classes of involutory automorphism of $G$ with $4$-dimensional 
centralizers (see \cite{lowen} (p. 155)). The centralizers of suitable representatives 
of these classes are isomorphic either to $H_1=Spin _3 \times SO_2(\mathbb R) / \langle (-1,-1) \rangle$ or to $H_2=SL_2(\mathbb R) \times SO_2(\mathbb R) / \langle (-1,-1) \rangle$.  Moreover, the $(-1)$-eigenspaces  ${\bf m}_i$, $i=1,2$ of these representatives are 

\medskip
\noindent
\centerline{${\bf m}_1=\langle e_4, e_5, e_7, e_8 \rangle$, \ \ 
${\bf m}_2=\langle e_1, e_3, e_4, e_5 \rangle$.} 

\medskip
\noindent
An Iwasawa decomposition (cf. \cite{iwasawa}, Theorem 6, p. 530) of the 
Lie algebra $\mathfrak{su_3}(\mathbb C,1)$ 
is given by 

\medskip
\noindent
\centerline{$\mathfrak{su_3}(\mathbb C,1)={\mathfrak k}+{\mathfrak a}+ 
{\mathfrak n}$,} 

\medskip
\noindent
where 
${\mathfrak k}=\langle e_1,e_2,e_3,e_6 \rangle $ is compact,  
${\mathfrak n}=\langle e_4-e_3, e_5+e_2, e_6+e_7 \rangle $ is nilpotent and 
${\mathfrak a}=\langle e_8 \rangle $. 
\newline
\noindent
Using this decomposition and the classification in \cite{chen}, Chap. 5, p. 276, we obtain that the conjugacy classes of the $4$-dimensional subalgebras of $\mathfrak{su_3}(\mathbb C,1)$ are the following 

\medskip
\noindent
${\bf h}_1=\mathfrak{so_3}(\mathbb R) \oplus \mathfrak{so_2}(\mathbb R)=\langle e_1, e_2, e_3, e_6 \rangle $, \  
${\bf h}_2=\mathfrak{sl_2}(\mathbb R) \oplus \mathfrak{so_2}(\mathbb R)=\langle e_2, e_6, e_7, e_8 \rangle $, 

\medskip
\noindent 
${\bf h}_3=\langle e_4-e_3, e_2+e_5, e_6+e_7, e_8 \rangle $.

\medskip
\noindent
The intersection of the subspace ${\bf m}_2$ with the subalgebras ${\bf h}_i$, $i=1,3$ as well as the intersection of ${\bf m}_1$ with the subalgebras 
${\bf h}_j$, $j=2,3$ is not trivial. Moreover, the subgroup $H_2$ can not be the stabilizer of the identity of a $4$-dimensional differentiable loop $L$ (see Lemma 5). Hence it remains to prove the triple $(G, H_1, \exp {\bf m}_1)$. 
\newline
The group $H_1$ is a maximal compact subgroup of $G$ hence there is a 
global 
differentiable Bruck loop $L_0$ of hyperbolic type  having 
$G \cong PSU_3(\mathbb C, 1)$ 
as the group topologically generated by its left translations. The  loop 
$L_0$ is realized on the 
complex  hyperbolic  plane geometry (cf. \cite{lowen}, p. 152). 
Since  only groups conjugate to $H_1$ can be complements of 
$\exp {\bf m}_1=\exp \{ \lambda _4 e_4+ \lambda _5 e_5+ \lambda _7 e_7+ 
\lambda _8 e_8;\ \lambda _i \in \mathbb R \}$ 
there is precisely one isotopism class 
${\cal C}$ of differentiable Bol loops  which are sections in 
$PSU_3(\mathbb C, 1)$. As a representative of ${\cal C}$ we can choose the 
loop $L_0$ which we call the complex  hyperbolic plane  loop.

\medskip
\noindent
The adjoint map  $\tau :X \mapsto (\bar{X})^t$ can be chosen as a representative of   involutions 
of ${\bf g}=\mathfrak{su_3}(\mathbb C,1)$ fixing elementwise  a  
$3$-dimensional subalgebra. The centralizer of 
$\tau $ is the subalgebra ${\bf h}=\langle e_2, e_4, e_7 \rangle 
\cong \mathfrak{sl_2}(\mathbb R)$.  The tangent space ${\bf m}$ of the 
$5$-dimensional symmetric space 
 corresponding to $\tau $ is generated by the basis elements 
$e_1, e_3, e_5, e_6, e_8$. 
\newline
Using the Iwasawa decomposition of $\mathfrak{su_3}(\mathbb C,1)$ and the classification in \cite{chen}, Chap. 5, p. 276, we see that every  $3$-dimensional  subalgebra 
${\bf h}$ of ${\bf g}$ has one of the following shapes:

\medskip
\noindent
\centerline{${\bf h_1} \cong \mathfrak{su_2}(\mathbb C)$, \ 
${\bf h_2} \cong \mathfrak{sl_2}(\mathbb R)$, }

\medskip
\noindent  
\centerline{${\bf h_3}=\langle e_5+e_2, e_6+e_7, e_8 \rangle $, \ 
${\bf h_4}=\langle e_4-e_3+b e_8, e_5+e_2, e_6+e_7 \rangle $, }

\medskip
\noindent  
\centerline{${\bf h_5}=\langle e_4-e_3+b (e_5+ e_2), e_6+e_7, e_8+c (e_5+e_2) 
\rangle $,}

\medskip
\noindent  
${\bf h_6}=\langle e_1-\frac{1}{2} e_6+\frac{3}{2} c(e_4-e_3)-\frac{3}{2} 
b(e_5+e_2), e_8+b(e_4-e_3)+c(e_5+e_2), 
\newline
e_6+e_7 \rangle$,  where $b, c \in \mathbb R$. 

\medskip
\noindent
Since the  subalgebras  of 
$\mathfrak{su_3}(\mathbb C,1)$ isomorphic to 
$\mathfrak{so_3}(\mathbb R)$ 
are conjugate under $Ad\ G$ (cf. \cite{chen}, Chap. 5, p. 276)
we may assume that 
${\bf h_1}=\langle e_1, e_2, e_3 \rangle $. 
The subalgebras of 
$\mathfrak{su_3}(\mathbb C,1)$ isomorphic to $\mathfrak{sl_2}(\mathbb R)$ form two conjugacy classes. As representatives of these conjugacy classes we can choose the subalgebras ${\bf h_{2,1}}=\langle e_2, e_4, e_7 \rangle $ and  
${\bf h_{2,2}}=\langle e_6, e_7, e_8 \rangle $.  
The  subalgebras ${\bf h_1}$ and  ${\bf h_{2,1}}$ 
 contain  the compact element $e_2$ which is conjugate to $e_1 \in {\bf m}$.
 The element 
$e_6+e_7 \in {\bf h_{2,2}} \cap {\bf h_3} \cap {\bf h_4}  \cap {\bf h_5} \cap {\bf h_6}$ is conjugate to 
$e_6+e_8 \in {\bf m}$ since  both are hyperbolic elements in  the same 
Lie subalgebras isomorphic to $\mathfrak{sl_2}(\mathbb R)$. 
\newline
The above  considerations yield the 

\begin{Theo} All Bol loops  having a group locally isomorphic to 
$PSU_3(\mathbb C,1)$  as the group 
topologically generated by its  left translations  are isotopic to the 
complex  hyperbolic plane  loop. 
\end{Theo}

\bigskip
\noindent
Let now $G$ be isomorphic to $SL_3(\mathbb R)$. According to 
\cite{lowen}, p. 155,  all involutory automorphisms of $SL_3(\mathbb R)$ are 
induced 
by reflections or polarities of the real projective plane  
${\mathbb P}_2(\mathbb R)$. 
The Lie algebra $\bf{g}=\mathfrak{sl_3}(\mathbb R)$ is isomorphic to the 
Lie 
algebra of  matrices
\[(\lambda _1 e_1+\lambda _2 e_2+\lambda _3 e_3+\lambda _4 e_4+\lambda _5 e_5+\lambda _6 e_6+\lambda _7 e_7+\lambda _8 e_8) \mapsto \]
\[ \left ( \begin{array}{ccc}
- \lambda _5-\lambda _8 & \lambda _1 & \lambda_2 \\
 \lambda _3 &  \lambda _5 & \lambda _6 \\
 \lambda _4 & \lambda _7 & \lambda _8 \end{array} \right ); \lambda_i \in 
\mathbb R, i=1,\cdots ,8. \] 
Then the Lie multiplication  of ${\bf g}$ is given by 
\[ [e_1,e_2]=[e_1,e_7]=[e_2,e_6]=[e_3,e_4]=[e_3,e_6]=[e_4,e_7]=[e_5,e_8]=0, \]
\[ [e_1,e_6]=[e_2,e_5]=\frac{1}{2} [e_2,e_8]=e_2,\  [e_1,e_8]=[e_2,e_7]=\frac{1}{2} [e_1,e_5]=e_1, \]
\[ [e_4,e_6]=[e_3,e_8]=\frac{1}{2} [e_3,e_5]=-e_3,\  [e_3,e_7]=[e_4,e_5]=\frac{1}{2} [e_4,e_8]=-e_4, \]
\[ [e_6,e_8]=[e_5,e_6]=[e_3,e_2]=e_6,\  [e_1,e_4]=[e_5,e_7]=[e_7,e_8]=-e_7,  \]
\[ [e_1,e_3]=- e_5,\  [e_2,e_4]=-e_8,\  [e_6,e_7]=e_5-e_8. \]
We choose an elliptic polarity in ${\mathbb P}_2(\mathbb R)$  such  
that 
the 
corresponding involution $\tau _1$ induces on $\mathfrak{sl_3}(\mathbb R)$ 
the  automorphism 
$\tau ^*_1: 
\mathfrak{sl_3}(\mathbb R) \to \mathfrak{sl_3}(\mathbb R);  X \mapsto -X^t$. 
The  Lie algebra 
$\mathfrak{so_3(\mathbb R)}= \langle e_1-e_3, e_2-e_4, e_7-e_6 \rangle$ is the $(+1)$-eigenspace of  
$\tau ^*_1$. 
The tangent space ${\bf m_1} \subset \mathfrak{sl_3}(\mathbb R)$  of the 
$5$-dimensional symmetric space 
belonging to $\tau _1^*$ has as generators $e_5$, $e_8$, $e_1+e_3$, 
$e_2+e_4$, $e_6+e_7$. We can choose a hyperbolic polarity in  
${\mathbb P}_2(\mathbb R)$  such  that the 
corresponding involution  $\tau _2^*$ of 
$\mathfrak{sl_3}(\mathbb R)$ 
is  given by $\tau _2^*:X \mapsto -\hbox{diag}(1,1,-1) X^t \ \hbox{diag}(1,1,-1)$. 
The Lie algebra fixed by $\tau _2^*$ elementwise is 
$\langle e_1-e_3, e_2+e_4, e_6+e_7 \rangle$. The tangent space ${\bf m_2}$ 
of the corresponding symmetric space   
is generated by $\{ e_1+e_3, e_2-e_4, e_6-e_7, e_5, e_8 \}$. 
\newline
According to Lemma 5 any maximal compact subgroup of the stabilizer $H$ is 
trivial or locally isomorphic to $SO_3(\mathbb R)$. 
Now using the classification  of Lie, who has determined all subalgebras of 
$\mathfrak{sl_3}(\mathbb R)$  (cf. \cite{lie}, pp. 288-289 and \cite{lie2}, 
p. 384) we obtain that  
the $3$-dimensional subalgebras ${\bf h}$ of the stabilizer $H$ have  
one of  the following shapes:

\medskip
\noindent
${\bf h_1}= \mathfrak{so_3}(\mathbb R)$, \ 
${\bf h_2}= \langle a(e_5+e_8)+e_6-e_7, e_1, e_2 \rangle $, $a >0$, \ 
${\bf h_3}= \langle e_1, e_2, e_6 \rangle $,

\medskip
\noindent
${\bf h_4}=  \langle e_5- e_8, e_2+ e_3, e_6, \rangle $,
${\bf h_5}=  \langle  e_3, e_6, e_8+e_2  \rangle $, 
${\bf h_6}= \langle  e_2, e_6,  e_5+ e_8- e_3 \rangle $, 

\medskip
\noindent
${\bf h_7}= \langle  e_5, e_8,  e_6 \rangle $, 
${\bf h_8}=  \langle  e_2, e_5+ e_8,  e_6  \rangle $, 
${\bf h_9}= \langle e_3, e_6, e_8 \rangle $,   

\medskip
\noindent
${\bf h_{10}}= \langle e_2, e_6, (b-1)e_5+b e_8, \rangle $, $b \in \mathbb R$,
${\bf h_{11}}= \langle e_3, e_6,  e_5+ c e_8 \rangle $, $c \in \mathbb R$.  

\medskip
\noindent
The element $e_8+e_2 \in {\bf h_{5}}$ is conjugate 
 to  $e_8 \in {\bf m_1} \cap {\bf m_2}$ under the element  
$g=\left( 
\begin{array}{rrr}
0 & 0 & 1 \\
0 & 1 & 0 \\
-1 & 0 & \frac{1}{2}  \end{array} \right)$. Futhermore the element  
$e_5+ e_8- e_3 \in {\bf h_6}$ is conjugate to $e_8-2 e_5 
\in {\bf m_1} \cap {\bf m_2}$ under   
$g=\left( 
\begin{array}{rrr}
0 & 0 & 1 \\
1 & 0 & 0 \\
-\frac{1}{3} & 1 & 0  \end{array} \right)$ 
and  $e_2 \in {\bf h_2} \cap {\bf h_{3}}$ is 
conjugate to 
$e_8-e_5-(e_6-e_7) \in {\bf m_2}$ under $g=\left( 
\begin{array}{rrr}
0 & -1 & 0 \\
1 & 0 & 0 \\
0 & 1 & 1  \end{array} \right) $.  
The intersection of the subspaces 
${\bf m_1}$  and ${\bf m_2}$ with the subalgebras ${\bf h_i}$, 
$i=4,7,8,9,10,11,$  as well as 
the intersection of ${\bf m_2}$ with the subalgebra ${ \bf h_1}$ is not 
trivial.  
Hence we may suppose that the Bol loop $L$ is realized on  
$\exp {\bf m_1}$ and 
the stabilizer of the 
identity of  $L$ has one of the following shapes: 

\medskip
\noindent
\centerline{a) \  $H_1=SO_3(\mathbb R)$,}

\medskip
\noindent
\centerline{b) \ $H_2=\left \{ \left ( 
\begin{array}{ccc} 
d^{-2 t} & a & b \\
0 & d^t \cos t & d^t \sin t \\
0 & -d^t \sin t & d^t \cos t \end{array} \right ); t \in [0, 2 \pi ), 
a,b \in \mathbb R \right \}$, with $d>1$,} 

\medskip
\noindent
\centerline{c) \  $H_3= \left \{ \left ( 
\begin{array}{ccc} 
1 & k  & l \\
0 & 1 & m \\
0 & 0 & 1 \end{array} \right ); k,l,m \in \mathbb R \right \}$. }

\begin{Prop} There is no differentiable Bol loop $L$ such that the 
stabilizer of the identity $e \in L$ is one of the Lie groups $H_i$, 
$i=2,3$ given in  b) and  c). 
\end{Prop}
\Pro The exponential image of the subspace ${\bf m_1}$ consists of all positive definite matrices of the shape 

\medskip
\noindent
\centerline{$\left \{ A= \left ( \begin{array}{ccc}
a & b & c \\
b & e & d \\
c & d & f \end{array} \right ); \det A=1 \right \}$. }

\medskip
\noindent
A differentiable Bol loop $L$ exists if and only if every coset 
$g H_i$, $g \in G$, $i=2,3$ contains precisely one element 
$m \in \exp {\bf m_1}$ (see  Lemma 6). For $c \in \mathbb R$ denote by 
$g_c H_i$ the coset 
$\left ( \begin{array}{ccc} 
1+c & 1 & 0 \\
c & 1 &  0 \\
0 & 0 & 1 \end{array} \right ) H_i$, $i=2,3$.  The coset $g_1 H_2$ 
contains different elements 
$m_1=\left ( \begin{array}{ccc} 
2 & 1 & 0 \\
1 & 1 &  0 \\
0 & 0 & 1 \end{array} \right )$ and  
$ m_2=\left ( \begin{array}{ccc} 
2 d^{-4 \pi} & d^{-4 \pi} & 0 \\
d^{-4 \pi} & \frac{ d^{-4 \pi} +d^{2 \pi }}{2}  &  0 \\
0 & 0 & d^{2 \pi} \end{array} \right )$ of $\exp {\bf m_1}$. 
If $c > -1$ then any coset $g_c H_3$  contains 
precisely the  element 

\medskip
\noindent
\centerline{$s_3(c)= \left ( \begin{array}{ccc} 
(1+c)  & c  & 0 \\
c   & \frac{c^2+1}{c+1}  &  0 \\
0 & 0 & 1 \end{array} \right ) \in \sigma (G/H_3)$.}  

\medskip
\noindent 
But $\lim\limits_{c \to -1} \sigma (g_c H_3)$ should  be $\lim\limits_{c \to -1} s_3(c)$ which is a contradiction.  
\qed

\bigskip
\noindent 
A reflection $\tau_3$  in $\mathbb P_2(\mathbb R)$ can be chosen in such a 
way that it  induces on 
$\mathfrak{sl_3}(\mathbb R)$ the involution 
$\tau_3^*: \mathfrak{sl_3}(\mathbb R) \to 
\mathfrak{sl_3}(\mathbb R); X \mapsto A^{-1} X A$,  
$A= \hbox{diag}\ (1,-1,-1)$, which  fixes  elementwise the 
subalgebra $\langle  e_5,e_6,e_7,e_8 \rangle$ isomorphic to 
 $\mathfrak{gl_2}(\mathbb R)$.
The tangent space ${\bf m}_3$ 
of the corresponding  $4$-dimensional symmetric space  has as 
generators $e_1$, $e_2$, $e_3$, $e_4$. 
The group associated with $\langle  e_5,e_6,e_7,e_8 \rangle$ is excluded by Lemma 5
and hence the classification  of Lie 
(cf. \cite{lie}, pp. 288-289 and \cite{lie2}, p. 384) yields that 
the Lie algebra ${\bf h}$ of the stabilizer of the identity of 
a $4$-dimensional Bol loop   has (up to conjugation) one of the following 
forms:  

\medskip
\noindent
${\bf h}_1=\langle e_1, e_2, e_6, e_5+c e_8 \rangle $, \   \  
${\bf h}_2=\langle e_3, e_5, e_6, e_8 \rangle $, \ \ 
${\bf h}_3=\langle e_1, e_2, e_6, e_8 \rangle $, 

\medskip
\noindent
${\bf h}_4=\langle e_2, e_5, e_6, e_8  \rangle $, 
where $c \in \mathbb R$. 

\medskip
\noindent
The intersection of all these subalgebras ${\bf h}_{i}$, $i=1, \cdots ,4$ 
with  the subspace ${\bf m}_3$ is not trivial.  

\medskip
\noindent
This contradiction to Lemma 6 and the above considerations yield the main 
part of the following

\begin{Theo} Every Bol loop with the group $SL_3(\mathbb R)$ as the group 
topologically generated by the  left translations is isotopic to the 
$5$-dimensional Bruck loop $L_0$ of hyperbolic type having the group 
$SO_3(\mathbb R)$ as the stabilizer  of $e \in L_0$.  
\end{Theo}       

\noindent
\Pro The group  $SO_3(\mathbb R)$ is a  maximal compact 
subgroup of 
$SL_3(\mathbb R)$. According to Theorem 8 there is a $5$-dimensional Bruck 
loop $L_0$ of hyperbolic type realized on the differentiable 
manifold $\exp {\bf m_1}$. 
Since up to isomorphisms there is only one symmetric space $\cal{S} $ having 
$SL_3(\mathbb R)$ as the group of displacements and  $SO_3(\mathbb R)$ as  
the centralizer of the  involutory 
automorphism belonging to $\cal{S} $
there exists  precisely one isotopism class 
${\cal C}$ of differentiable Bol loops corresponding to  $SL_3(\mathbb R)$  
and $L_0$ is a representative of ${\cal C}$. 
\qed

\section{ Bol loops corresponding to semi-simple Lie groups} 
\label{sec:4}

\bigskip
\noindent
Let $G=G_1 \times G_2$  be the group topologically generated by the left 
translations of a 
connected differentiable Bol loop $L$ such that $M=M_1 \times M_2$ holds, where $M$ is the symmetric space belonging to $L$. 
If for the stabilizer $H$ of 
$e \in L$ one has $H=H_1 \times H_2$ with  $1 \neq H_i <  G_i$, $i=1,2,$ then $L$ is a direct product of two proper Bol loops $L_1$ and $L_2$
such that $L_i$ 
is  realized 
on $M_i$,  has  $G_i$ as the group topologically generated by the left 
translations and $H_i$ as the stabilizer of $e \in L_i$, $i=1,2,$ 
(cf. \cite{loops}, 
Proposition 1.19, p.28).  If $M_2 \cong G_2$ and the stabilizer $H$ of $e \in L$ has the shape $H=(H_1, \varphi (H_1))$, where $\varphi :H_1 \to G_2$ is a homomorphism,  then $L$ is a Scheerer extension of the Lie group $G_2$ by the proper Bol loop $L_1$  (cf. \cite{loops}, Proposition 2.4, p. 44).  
If $G$ has dimension $\le 9$ and the direct factors are simple  then $G_i$ 
is isomorphic either to $PSL_2(\mathbb R)$ or to $PSL_2(\mathbb C)$ or to $SO_3(\mathbb R)$.  There 
is no proper Bol loop corresponding to $SO_3(\mathbb R)$ (cf. \cite{loops}, 
Corollary 16.8, p. 204) and  every proper Bol loop having $PSL_2(\mathbb R)$ 
respectively $PSL_2(\mathbb C)$ as the group topologically generated by the 
left translations is isotopic to the hyperbolic plane loop 
(cf. \cite{loops}, Section 22) respectively to the hyperbolic space loop
(cf. \cite{figula1}, Theorem 5, p. 446 and Proposition 10).
 This discussion yields 

\begin{Theo} Let $L$ be a connected differentiable Bol loop having an at 
most $9$-dimensional semi-simple Lie group $G$ as the group topologically 
generated by its left translations. 
\newline
i) If the stabilizer $H$ is a direct product of subgroups $1 \neq H_i$ contained in the simple factors $G_i$ of $G$ then $L$ is a direct product of 
proper Bol loops $L_i$ isotopic to the hyperbolic plane loop $\mathbb H_2$ respectively to the hyperbolic space loop $\mathbb H_3$. Furthermore, 
$G$ is isomorphic either to  
$PSL_2(\mathbb R) \times PSL_2(\mathbb R)$ or to  
$PSL_2(\mathbb R) \times PSL_2(\mathbb R) \times PSL_2(\mathbb R)$ 
respectively to $PSL_2(\mathbb C) \times PSL_2(\mathbb R)$. 

\smallskip
\noindent
ii) If the stabilizer $H$ is not a direct product of subgroups $1 \neq H_i$ contained in the simple factors $G_i$ of $G$ then one has 
$G=G_1 \times S$, where $G_1$ is isomorphic either to $PSL_2(\mathbb R)$ or to $PSL_2(\mathbb R) \times PSL_2(\mathbb R)$ respectively to $PSL_2(\mathbb C)$ and $S$ is the complement of $G_1$ in $G$. Moreover, the stabilizer $H$ has the shape $\{(x, \varphi (x))\ |\ x \in H_1 \}$, 
where $H_1$ is isomorphic either to $SO_2(\mathbb R)$ or to $SO_2(\mathbb R) \times SO_2(\mathbb R)$ 
respectively to $SO_3(\mathbb R)$ and  the loop $L$ is 
 a Scheerer extension of $S$ by $\mathbb H_2$ in the first case, by $\mathbb H_2 \times \mathbb H_2$ in the second case and by $\mathbb H_3$ 
 in the third case.  
\end{Theo}

\noindent From now we assume that the subgroup $H$ is not decomposable 
into a direct product. Denote by $p_i:G \to G_i$ the projection of $G$ onto 
the $i$-th components $G_i$ of $G$. 

\begin{Lemma} There is one conjugacy class ${\mathcal C_1}$ of  involutory 
automorphisms of $\mathfrak{so_3}(\mathbb R)$ and two conjugacy classes 
${\mathcal C_2}$ and ${\mathcal C_3}$ of 
involutory automorphisms of $\mathfrak{sl_2}(\mathbb R)$. As a 
representative of ${\mathcal C_1}$  we can choose  one which 
fixes the $1$-dimensional subalgebra $\langle e_3 \rangle$  elementwise. 
As a 
representative of ${\mathcal C_2}$ respectively ${\mathcal C_3}$ we can 
choose  one  
fixing  $\langle e_3 \rangle$ respectively 
$\langle e_2 \rangle$ elementwise. 
\end{Lemma}

\noindent
\Pro The assertion follows from \cite{furness}, pp. 44-45. 
\qed

\begin{Prop} Every  connected differentiable Bol loop having a group locally 
isomorphic to 
$PSL_2(\mathbb R) \times SO_3(\mathbb R)$ as the group topologically 
generated by its left translations is isotopic to a Scheerer extension of $SO_3(\mathbb R)$ by $\mathbb H_2$.   
\end{Prop} 

\noindent
\Pro Assume that $L$ is not a Scheerer extension in the assertion.    
The automorphism group $\Gamma $ of the Lie algebra ${\bf g}$ of $G$ is 
the 
direct product of the automorphism group 
of $\mathfrak{sl_2}(\mathbb R)$ and the 
automorphism group of $\mathfrak{so_3}(\mathbb R)$. 
According to Lemma 15  there are two conjugacy classes of involutory 
automorphisms of ${\bf g}$ fixing elementwise a $2$-dimensional subalgebra.  
The $(-1)$-eigenspaces of suitable representatives of these classes are 

\medskip
\noindent
${\bf m}_1=\langle (e_1,0), (e_2,0), (0,\hbox{i} e_1), (0,\hbox{i} e_2) 
\rangle$, 
${\bf m}_2=\langle (e_1,0), (e_3,0), (0,\hbox{i} e_1), (0,\hbox{i} e_2) 
\rangle $.  

\medskip
\noindent  
Since $SO_3(\mathbb R)$ has no $2$-dimensional subgroup  we have 
to investigate the case that  $\hbox{dim}\  p_1(H) = 2$ 
and $\hbox{dim}\  p_2(H) = 1$. Then  we may assume that  
$p_1({\bf h})=\langle  e_2+e_3, e_1 \rangle $, 
$p_2({\bf h})=\langle e_3 \rangle$ (see  {\bf 1.1} and  {\bf 1.2}) 
and  the Lie algebra  
 ${\bf h}$ has the shape ${\bf h}=\langle (e_2+e_3,0), (e_1,e_3) \rangle $. 
The element $(e_2+e_3,0) \in {\bf h}$ is conjugate to 
$(e_1+e_3,0) \in {\bf m}_2$ (see {\bf 1.1}). Hence there is no 
differentiable Bol loop $L$ with $T_e L={\bf m}_2$. 
Since $p_1(\exp {\bf m}_1)$ and $p_1(\exp {\bf h})$ satisfy the conditions 
of Lemma 7 we have also here a contradiction. 
\qed

\bigskip
\noindent
\begin{Prop} If  the Lie 
group $G'=G_1 \times G_2$ or  $G''=G_1 \times G_2 \times G_3$, where  $G_i$ 
$(i=1,2,3)$ is  locally isomorphic to 
$PSL_2(\mathbb R)$, is  the group topologically generated by the  left 
translations of a connected differentiable Bol loop $L$ then $L$ is either a direct product  of proper Bol loops isotopic to the hyperbolic plane loop $\mathbb H_2$ or a Scheerer extension of $PSL_2(\mathbb R)$ by $\mathbb H_2$ respectively by $\mathbb H_2 \times \mathbb H_2$ or a Scheerer extension of $PSL_2(\mathbb R) \times PSL_2(\mathbb R)$ by $\mathbb H_2$. 
\end{Prop}
\Pro Assume that $L$ has not a form as in the assertion.   
The automorphism group 
$\Gamma $ of the Lie algebra ${\bf g'}=\mathfrak{sl_2}(\mathbb R) \oplus 
\mathfrak{sl_2}(\mathbb R)$ of $G'$ is the 
semidirect product  of the normal automorphism group $\Gamma _1 
\times  \Gamma_1 $, where $\Gamma _1$ is the automorphism group of 
$\mathfrak{sl_2}(\mathbb R)$, 
by the group generated by  
$\sigma :\mathfrak{sl_2}(\mathbb R) \to \mathfrak{sl_2}(\mathbb R); 
\ (u,v) \mapsto (v,u)$. To the symmetric space determined by $\sigma $ 
there corresponds a $3$-dimensional Bol loop having 
$PSL_2(\mathbb R) \times PSL_2(\mathbb R)$ as the group topologically 
generated by the left translations, but such a Bol loop does not exit  
(cf. \cite{figula1}, pp. 442-444). 
 Hence  there exist up to  
 automorphisms of $G'$ precisely three $4$-dimensional  symmetric 
spaces of $G'$  the tangent spaces  of which are given by  

\medskip
\noindent
\centerline{${\bf m}_1=\langle (e_1,0), (e_2,0), (0,e_1), (0,e_2) \rangle $, 
${\bf m}_2=\langle (e_1,0), (e_2,0), (0,e_1), (0,e_3) \rangle $,}

\smallskip
\noindent
\centerline{${\bf m}_3=\langle (e_1,0), (e_3,0), (0,e_1), (0,e_3) \rangle $.}

\medskip
\noindent
Moreover,  there are  up to automorphisms of 
$G''$   four $6$-dimensional symmetric 
spaces in $G''$  the tangent spaces of which are  

\medskip
\noindent
\centerline{${\widetilde {\bf m}}_1=\langle (e_1,0,0), (e_2,0,0), (0,e_1,0), (0,e_2,0), (0,0,e_1), (0,0,e_2) \rangle $, }

\smallskip
\noindent
\centerline{${\widetilde {\bf m}}_2=\langle (e_1,0,0), (e_2,0,0), (0,e_1,0), (0,e_2,0), 
(0,0,e_1), (0,0,e_3) \rangle $, }

\smallskip
\noindent
\centerline{${\widetilde {\bf m}}_3=\langle (e_1,0,0), (e_2,0,0), (0,e_1,0), (0,e_3,0), (0,0,e_1), (0,0,e_3) \rangle $, }

\smallskip
\noindent
\centerline{${\widetilde {\bf m}}_4=\langle (e_1,0,0), (e_3,0,0), (0,e_1,0), (0,e_3,0), (0,0,e_1), (0,0,e_3) \rangle $.}

\medskip 
\noindent
If $L$ corresponds to $G'$ then $\hbox{dim}\ H=2$ and 
$H=(\varphi({\mathcal L}_2), {\mathcal L}_2)$, where  
$\varphi \neq 1$ is a homomorphism of ${\mathcal L}_2$
into  $PSL_2(\mathbb R)$. If $\varphi $ is injective  then the Lie algebra 
of $H$ has the shape 
${\bf h}= \langle (e_1,e_1), (e_2+e_3, e_2+e_3) \rangle $ and the 
intersection of ${\bf h}$ with  ${\bf m}_i$ 
for  $i=1,2,3,$ is not trivial. 
If $\varphi $ has $1$-dimensional kernel  then ${\bf h}$ contains the 
element $(0,e_2+e_3)$ which is conjugate to   
 $(0,e_1+e_3) \in {\bf m}_2 \cap {\bf m}_3$ (see  {\bf 1.1}).  
Hence the subspaces  ${\bf m}_2$ and 
${\bf m}_3$ cannot determine  a Bol loop. Since 
$p_2(\exp {\bf m}_1)$ and $p_2(\exp {\bf h})$ have a  shape as in Lemma 7 
we obtain also a contradiction. 

\smallskip
\noindent
If $L$ corresponds to $G''$ then $\hbox{dim}\ H \in \{ 3,4 \}$ 
(see Lemma 1).  
The dimension of  $H$ cannot be $4$ since every 
$4$-dimensional subgroup of $G''$ contains a 
direct factor and  there is no $3$-dimensional Bol loop corresponding to 
$PSL_2(\mathbb R) \times PSL_2(\mathbb R)$.

\smallskip
\noindent
Let now $\hbox{dim}\ H=3$ and  
$p_i(H) \neq (H \cap G_i)$ for  $i=1,2,3$. 
If $p_1(H)$ is 
isomorphic to $PSL_2(\mathbb R)$ then one has 
$p_j(H) \cong PSL_2(\mathbb R)$ for $j=2,3$ and  the Lie algebra 
${\bf h}$ of $H$ contains  $(e_1,e_1,e_1)$ up to conjugacy. If $p_i(H)$ is 
isomorphic to 
${\mathcal L}_2$  for  $i=2,3$  
then the projection of 
$H$ onto the first  component of $G''$ is a non-trivial  homomorphic image 
of ${\mathcal L}_2$. In this case the Lie algebra 
${\bf h}$ contains either  the element $(e_1,e_1,e_1)$ or 
$(0,0,e_2+e_3)$ up to conjugacy. The first of them  lies in 
$\widetilde{{\bf m}_i}$ for  $i=1,2,3,4$ and the second is conjugate to 
$(0,0,e_1+e_3) \in \widetilde{{\bf m}_i}$, $i=2,3,4$ (see {\bf 1.1}). 
Since $p_3(\exp \widetilde{{\bf m}_1})$ and $p_3(\exp {\bf h})$ satisfies  
the conditions of Lemma 7 we have a contradiction and the assertion follows. 
\qed

\bigskip
\noindent
\begin{Prop} Every connected differentiable  Bol loop $L$ having a group  
locally isomorphic either to 
$G'=PSL_2(\mathbb R) \times SO_3(\mathbb R) \times SO_3(\mathbb R)$ or to 
$G''=PSL_2(\mathbb R) \times PSL_2(\mathbb R) \times SO_3(\mathbb R)$ as 
the group topologically generated by its left translations is a Scheerer extension of $SO_3(\mathbb R) \times SO_3(\mathbb R)$ respectively of 
$PSL_2(\mathbb R) \times SO_3(\mathbb R)$ by the hyperbolic plane loop $\mathbb H_2$ or a Scheerer extension of $SO_3(\mathbb R)$ by $\mathbb H_2 \times \mathbb H_2$. 
\end{Prop}

\noindent
\Pro Assume that $L$ is not a Scheerer extension in the assertion. 
Lemmata 1, 2 and 3  exclude the group $G'$. 
\newline
According to Proposition 16 and Lemma 2 we may assume that  
$p_i(H) \neq (H \cap G_i)$, 
$i=1,2,3$. Hence the stabilizer $H$ of the identity of 
$L$ in  $G''$ has dimension $3$. Moreover,  the Lie algebra of $H$ has 
the form 

\medskip
\noindent
\centerline{${\bf h}= \langle (e_1,e_1,e_3), (e_2+e_3,0,0), (0,e_2+e_3,0) 
\rangle.$ }

\medskip
\noindent
According to Lemma 15  we have four conjugacy classes of 
involutory automorphisms of 
${\bf g}''=sl_2(\mathbb R) \oplus sl_2(\mathbb R) \oplus so_3(\mathbb R)$ 
which fix elementwise a $3$-dimensional subalgebra of ${\bf g}''$. 
The $(-1)$-eigenspaces of suitable representatives of these classes are 
given by 

\medskip
\noindent
\centerline{${\bf m}_1=\langle (e_1,0,0), (e_2,0,0), (0,e_1,0), (0,e_2,0), 
(0,0,\hbox{i} e_1), (0,0,\hbox{i} e_2) \rangle$, }

\medskip
\noindent
\centerline{${\bf m}_2=\langle (e_1,0,0), (e_2,0,0), (0,e_1,0), (0,e_3,0), 
(0,0,\hbox{i} e_1), (0,0,\hbox{i} e_2) \rangle$, }

\medskip
\noindent
\centerline{${\bf m}_3=\langle (e_1,0,0), (e_3,0,0), (0,e_1,0), (0,e_2,0), 
(0,0,\hbox{i} e_1), (0,0,\hbox{i} e_2) \rangle$, }

\medskip
\noindent
\centerline{${\bf m}_4=\langle (e_1,0,0), (e_3,0,0), (0,e_1,0), (0,e_3,0), 
(0,0,\hbox{i} e_1), (0,0,\hbox{i} e_2) \rangle$. }

\medskip
\noindent
As the elements  $e_3$ and $\hbox{i} e_1$ are conjugate in 
$SO_3(\mathbb R)$ (see  {\bf 1.2}) the subalgebras ${\bf m}_i$, $i=1,2,3,4,$ 
cannot be the tangent spaces of  Bol loops  corresponding to $G''$. 
\qed

\begin{Prop} Let the group $G$ be locally isomorphic to a direct product 
$G=G_1 \times G_2$, where $G_1 \cong PSL_2(\mathbb C)$ and $G_2$ is either 
$PSL_2(\mathbb R)$ or $SO_3(\mathbb R)$.  
Every connected differentiable Bol loop  having $G$ as the 
group topologically generated by its left translations is either a direct product 
of the hyperbolic space loop $\mathbb H_3$ and the hyperbolic plane loop $\mathbb H_2$ or a Scheerer extension of $PSL_2(\mathbb C)$ by $\mathbb H_2$ or a Scheerer extension of a $3$-dimensional simple Lie group by $\mathbb H_3$.  
\end{Prop} 
\Pro 
There is precisely one conjugacy class of $4$-dimensional subgroups of 
$SL_2(\mathbb C)$ 
(see \cite{asoh}, p. 277). This class can be represented by a  subgroup 
$V$, the Lie algebra of which is  

\medskip
\noindent
\centerline{$v=\langle e_1, \hbox{i}\ e_1,e_2+e_3, \hbox{i}(e_2+e_3) 
\rangle$.}

\medskip
\noindent
The maximal compact subgroups  of $V$ are isomorphic to $SO_2(\mathbb R)$. 
Representatives of the $3$-dimensional subalgebras of $SL_2(\mathbb C)$ (cf. 
\cite{asoh}, p. 277-278) are given by

\medskip
\noindent
\centerline{$\mathfrak{so_3}(\mathbb R)= \langle \hbox{i} e_1, \hbox{i} e_2, e_3 
\rangle$, \ \  $\mathfrak{sl_2}(\mathbb R)= \langle e_1, e_2, e_3 
\rangle$,}

\medskip
\noindent
\centerline{$w_r= \langle  (r \hbox{i} -1)e_1, e_2+e_3, 
\hbox{i}(e_2+e_3) \rangle $, \ \ $u_1= \langle \hbox{i} e_1, e_2+e_3, 
\hbox{i}(e_2+e_3) \rangle $.}

\medskip
\noindent
Any $2$-dimensional ideal of the subalgebras $w_r$, $r \in \mathbb R$, and 
$u_1$ is isomorphic to $\langle e_2+e_3, \hbox{i} (e_2+e_3) \rangle $. 
The maximal compact subgroups of the Lie group  corresponding to $w_r$ are 
trivial for any $r \in \mathbb R$, whereas the maximal compact subgroups of 
the Lie group corresponding to $u_1$ are isomorphic to $SO_2(\mathbb R)$. 
\newline
According to Lemma 1  we have $ \hbox{dim}\ L \in \{5,6 \}$.   
If $G$ is locally isomorphic to $PSL_2(\mathbb C) \times SO_3(\mathbb R)$,  
then $ \hbox{dim}\ L \neq 5$. Otherwise $\hbox{dim}\ H$ would be $4$ 
and $p_1(H)=V$  which contradicts Lemma 5. 
\newline
There are precisely three classes of involutions of 
$\mathfrak{sl_2}(\mathbb C)$ 
(cf. \cite{lowen}, pp. 152-153). The $(+1)$-eigenspaces of suitable 
representatives of these involutions are  
${\bf h}=\langle e_3, i e_3 \rangle$ or  
${\bf h}=\langle e_1, e_2, e_3 \rangle $ 
respectively ${\bf h}=\langle \hbox{i} e_1, \hbox{i} e_2, e_3 \rangle $.  
The automorphism group $\Gamma $ of the Lie algebra 
${\bf g}'=\mathfrak{sl_2}(\mathbb C) \oplus \mathfrak{sl_2}(\mathbb R)$ 
respectively  of 
${\bf g}''=\mathfrak{sl_2}(\mathbb C) \oplus \mathfrak{so_3}(\mathbb R)$
is the direct product of the 
automorphism group of $\mathfrak{sl_2}(\mathbb C)$ and  the automorphism 
group of $\mathfrak{sl_2}(\mathbb R)$ respectively 
of $\mathfrak{so_3}(\mathbb R)$.  
According to Lemma 15  there exist precisely two conjugacy classes of 
involutions of ${\bf g}'$ having  $3$-dimensional subalgebras  as their 
$(+1)$-eigenspaces. The $(-1)$-eigenspaces of suitable representatives of 
these classes are  given by

\medskip
\noindent
\centerline{${\bf m}_1=\langle (e_1,0), (e_2,0), (\hbox{i} e_1,0), 
(\hbox{i} e_2,0), (0,e_1), (0,e_2) \rangle$, }  

\medskip
\noindent
\centerline{${\bf m}_2=\langle (e_1,0), (e_2,0), (\hbox{i} e_1,0), 
(\hbox{i} e_2,0), (0,e_1), (0,e_3) \rangle$. } 

\medskip
\noindent
Moreover, there are  precisely four conjugacy classes of 
involutions of ${\bf g}'$ leaving  $4$-dimensional subalgebras  elementwise 
fix. The tangent spaces ${\bf m}_i$, $i=3,4,5,6,$ of the corresponding  
symmetric spaces are:  

\medskip
\noindent
\centerline{${\bf m}_3=\langle (\hbox{i} e_1,0), (\hbox{i} e_2,0), 
(\hbox{i} e_3,0), (0,e_1), (0,e_2) \rangle,$ } 

\medskip
\noindent
\centerline{${\bf m}_4=\langle (\hbox{i} e_1,0), (\hbox{i} e_2,0), 
(\hbox{i} e_3,0), (0,e_1), (0,e_3) \rangle$, } 

\medskip
\noindent
\centerline{${\bf m}_5=\langle (e_1,0), (e_2,0), (\hbox{i} e_3,0), (0,e_1), 
(0,e_2) \rangle$, } 

\medskip
\noindent
\centerline{${\bf m}_6=\langle (e_1,0), (e_2,0), (\hbox{i} e_3,0), (0,e_1), 
(0,e_3) \rangle$. } 

\medskip
\noindent
For ${\bf g}''$  there is one conjugacy class of involutory automorphisms 
 which fix   $3$-dimensional subalgebras  elementwise. The 
$(-1)$-eigenspace of a  representative  of this class is the subspace 

\medskip
\noindent
\centerline{${\widetilde {\bf m}}_1=\langle (e_1,0), (e_2,0), 
(\hbox{i} e_1,0), (\hbox{i} e_2,0), (0,\hbox{i} e_1), (0,\hbox{i} e_2) 
\rangle$. }

\medskip
\noindent 
First we consider the case that $G$ is locally isomorphic to 
$PSL_2(\mathbb C) \times SO_3(\mathbb R)$ and  $\hbox{dim}\ H=3$. According 
to Lemma 5 the 
projection $p_1({\bf h})$ of 
the Lie algebra ${\bf h}$ of $H$ onto the first component is isomorphic 
either to 
$\mathfrak{so_3}(\mathbb R)$ or to a Lie algebra 
$w_r$,  $r \in \mathbb R$. 
If $p_1({\bf h}) \cong p_2({\bf h}) \cong \mathfrak{so_3}(\mathbb R)$ then 
we may assume that ${\bf h}$ has the shape $\{ (x,x);\ x \in 
\mathfrak{so_3}(\mathbb R) \}$.  But then  
${\bf h} \cap {\widetilde {\bf m}}_1$ is  not trivial. If 
$p_1({\bf h})$ is isomorphic 
to $w_r$, $r \in \mathbb R$, then up to conjugation ${\bf h}$ has the form 
${\bf h}=(w_r, \phi (w_r))$, where $\phi $ is a homomorphism of $w_r$ onto 
$\mathfrak{so_2}(\mathbb R)$  and 
$\phi ^{-1}(0)=\langle (e_2+e_3,0),  
(\hbox{i} (e_2+e_3,0)) \rangle $. Since $(e_2+e_3,0)$ is conjugate to 
$(e_1+\hbox{i} e_2,0) \in {\widetilde {\bf m}}_1$ the subspace 
${\widetilde {\bf m}}_1$ cannot be the tangent space of a differentiable 
Bol loop. 
\newline
Now we assume that $G$ is locally isomorphic to 
$PSL_2(\mathbb C) \times PSL_2(\mathbb R)$ and  $\hbox{dim}\ H=3$.  
If $\hbox{dim}\ p_2({\bf h})=3$ then the subgroup $H$ is locally isomorphic 
to  
$\{ (x,x); \ x \in PSL_2(\mathbb R) \}$ and the Lie algebra ${\bf h}$ of $H$ is generated by 
$(e_1,e_1)$, $(e_2+e_3,e_2+e_3)$, $(e_3,e_3)$. But then     
${\bf h} \cap {\bf m}_i \neq \{ 0 \}$ for $i=1,2$. 
\newline
If $\hbox{dim}\ p_2({\bf h})=2$ then 
$\hbox{dim}\ p_1({\bf h}) \in \{2,3 \}$. If 
$ p_1(H) \cong p_2(H) \cong \mathcal{L}_2$ then the Lie algebra ${\bf h}$ 
has the shape $\langle (e_1,e_1),(e_2+e_3,0), (0,e_2+e_3) \rangle$. If  
$p_1({\bf h})$ is a $2$-dimensional abelian subalgebra of 
$\mathfrak{sl_2}(\mathbb C)$ then 
 ${\bf h}$ has the shape 
${\bf h}=(K,0) \oplus (\varphi (p_2({\bf h})), p_2({\bf h}))$, where 
$\varphi $ is a 
homomorphism with the nucleus $(0,e_2+e_3)$ and $K$ is a complement of 
$\varphi (p_2({\bf h}))$ in $p_1({\bf h})$. According to Lemma 9 the 
homomorphism $\varphi $ may be chosen 
in such a way that $\varphi (p_2({\bf h}))$ has one of the following shapes: 
$\langle e_3 \rangle$, $\langle \hbox{i} e_3 \rangle$, 
$\langle e_2+e_3 \rangle$, $\langle \hbox{i} (e_2+e_3) \rangle$. Hence  
${\bf h}$ contains one of the following  $1$-dimensional 
algebras: $\langle (e_3,e_1) \rangle $, 
$\langle (\hbox{i} e_3,e_1) \rangle $, 
$\langle (e_2+e_3,e_1) \rangle $, 
$\langle ( \hbox{i} (e_2+e_3),e_1) \rangle $.    
But $(e_1,e_1) \in {\bf h} \cap {\bf m}_1 \cap {\bf m}_2$ and 
  the  element $(0,e_2+e_3) \in {\bf h} $ is conjugate to 
$(0,e_1+e_3) \in {\bf m}_2$ (see {\bf 1.1}). Since $p_1({\bf m}_1)$ as well 
as $p_1({\bf h})$ have a  form as in  Lemma 9 these subalgebras ${\bf h}$ 
are excluded. 
\newline
Let now  $\hbox{dim}\  p_1({\bf h})=3$ and 
$\hbox{dim}\  p_2({\bf h}) \in \{2,1 \}$.  Then ${\bf h}$ has the shape 
${\bf h}=(p_1({\bf h}), \varphi (p_1({\bf h})))$, where 
$\varphi $ is a homomorphism of  a $3$-dimensional subalgebra of 
$\mathfrak{sl_2}(\mathbb C)$ onto $p_2({\bf h})$ with 
$\hbox{dim}\ \varphi ^{-1}(0) \in \{2,1 \}$. Then $p_2({\bf h})$ 
contains the element $(e_2+e_3,0)$ or $(\hbox{i} (e_2+e_3,0))$. However, 
both 
of them  are conjugate 
to $(e_1+\hbox{i} e_2,0) \in {\bf m}_1 \cap {\bf m}_2$ (see Lemma 9). 

\medskip
\noindent
Finally let $G$ be locally isomorphic to 
$PSL_2(\mathbb C) \times PSL_2(\mathbb R)$ and let $\hbox{dim}\ H=4$. Since 
$H$ does not decompose onto 
a direct product $H_1 \times H_2$ with  $H_i < G_i$ we have 
$\hbox{dim}\  p_1({\bf h})+ \hbox{dim}\  p_2({\bf h}) \ge 5$ and  hence   
$\hbox{dim}\  p_1({\bf h}) \ge 3$. 
\newline
We consider 
 first the case $\hbox{dim}\  p_1({\bf h})= 3$. The subalgebra 
$p_2({\bf h})$ cannot be $\mathfrak{sl_2}(\mathbb R)$  since 
$\mathfrak{sl_2}(\mathbb R)$ is simple. 
\newline
If  
$\hbox{dim}\  p_2({\bf h})=2$ then  we may assume that $p_2({\bf h})=
\langle e_1, e_2+e_3 \rangle $ (see {\bf 1.1}). Then  ${\bf h}$ has the 
form  
$(p_1({\bf h}), \varphi (p_1({\bf h})) \oplus \langle (0,e_2+e_3) \rangle$, 
where $\varphi $ is a 
homomorphism of  $p_1({\bf h})$ into $p_2({\bf h})$ with 
$\hbox{dim} \  \varphi ^{-1}(0)=2$ and  
 $p_1({\bf h})$ is either the subalgebra 
$u_1$ or  $w_r$, $r \in \mathbb R$. It follows that  ${\bf h}$ 
contains the elements $(e_2+e_3,0)$ and  $(\hbox{i} (e_2+e_3),0)$. The 
element $(\hbox{i} (e_2+e_3),0)$ is contained in  ${\bf m}_i$, $i=3,4$, 
and  $(0,e_1+e_3) \in {\bf m}_6$ is conjugate to  
$(0,e_2+e_3) \in {\bf h}$ (see  {\bf 1.1}).  
\newline
Therefore it remains to 
investigate the triples 
$(G,\exp {\bf h}_i, \exp {\bf m}_5)$, $i=1,2$, with  

\medskip
\noindent
\centerline{${\bf h}_1=\langle (e_2+e_3,0), (\hbox{i} (e_2+e_3),0), ((r \hbox{i}-1)e_1,e_1), (0,e_2+e_3) \rangle $}

\medskip
\noindent
\centerline{${\bf h}_2=\langle (e_2+e_3,0), (\hbox{i} (e_2+e_3),0), 
(\hbox{i}  e_1,e_1), 
(0,e_2+e_3) \rangle $. }

\medskip
\noindent  
Let  $\hbox{dim} \ p_1({\bf h})= 4$. Then up to conjugation 
we have $p_1({\bf h})=v$. 
Since $v$ is solvable and $\mathfrak{sl_2}(\mathbb R)$ is simple 
it follows $\hbox{dim}\ p_2({\bf h})<3$. 
\newline
If  $\hbox{dim} \ p_2({\bf h})=2$ then we may assume that 
$p_2({\bf h})=\langle e_1, e_2+e_3 \rangle $ (see {\bf 1.1}) and the 
subalgebra ${\bf h}$ has the form 
${\bf h}=(v, \psi (v))$, where $\psi $ is a homomorphism from $v$ onto  
$p_2({\bf h})$ 
such that $\psi ^{-1}(0)=\langle e_2+e_3, \hbox{i}(e_2+e_3) \rangle$. Since 
the image $\psi (w_r)$, 
$r \in \mathbb R$, as well as  $\psi (u_1)$ is the $1$-dimensional ideal 
$\langle e_2+e_3 \rangle $ of 
$p_2({\bf h})$ the subalgebra ${\bf h}$ would  contain  
$((r \hbox{i}-1) e_1,e_2+e_3)$ for all $r \in \mathbb R$ and $(\hbox{i} e_1,e_2+e_3)$. 
This contradicts  $\hbox{dim}\ {\bf h}=4$.  
\newline
Let now $\hbox{dim} \ p_2({\bf h})=1$. Then one has 
${\bf h}=(v, \varphi (v))$, where $\varphi $ is a homomorphism from $v$ onto 
a $1$-dimensional subalgebras of $\mathfrak{sl_2}(\mathbb R)$.  
Since 
$ \varphi^{-1}(0)$ is $3$-dimensional the  commutator subalgebra 
$\langle (e_2+e_3,0), (\hbox{i}(e_2+e_3),0) \rangle$ of $v$ is contained in 
$ \varphi^{-1}(0)$.  
If the element $((r \hbox{i}-1) e_1,0) \in w_r$ lies in ${\bf h}$ then for  
the 
fourth generator of ${\bf h}$ we have up to conjugation the following 
possibilities:
$(e_1,e_1)$, $(e_1,e_2+e_3)$, $(e_1,e_3)$, $(\hbox{i} e_1,e_1)$, 
$(\hbox{i} e_1,e_2+e_3)$, 
$(\hbox{i} e_1,e_3)$. If ${\bf h}$ 
contains  the element $(\hbox{i} e_1,0)$  
then for the fourth basis element of ${\bf h}$ we may choose 
 one of the following:   
$(e_1,e_1)$, 
$(e_1,e_2+e_3)$, $(e_1,e_3)$.  The element 
$(\hbox{i}(e_2+e_3),0) \in {\bf h}$ lies  in  
${\bf m}_3 \cap {\bf m}_4$,  the 
elements 
$(e_1,e_1)$ and $((r \hbox{i}-1) e_1,0)-r(\hbox{i} e_1, e_1)$ of ${\bf h}$ 
are 
contained in ${\bf m}_5 \cap {\bf m}_6$ and  $(e_1,e_3)$, 
$((r \hbox{i} -1) e_1,0)-r( \hbox{i} e_1, e_3)$ are elements of ${\bf m}_6$. Moreover, 
$(e_1,e_2+e_3)$ respectively  
$((r \hbox{i}-1) e_1,0)-r( \hbox{i} e_1,e_2+e_3)$ is  conjugate to  
$(e_1,e_1+e_3)$ 
respectively to 
$(-e_1,-r(e_1+e_3))$ which are elements of  ${\bf m}_6$ (see {\bf 1.1}). 
\newline
It 
remains to 
investigate the 
triples 
$(G, \exp {\bf h}_i, \exp {\bf m}_5)$, where 
${\bf h}_i$ has one of the following shapes: 

\medskip
\noindent
${\bf h}_3= \langle (e_2+e_3,0), (\hbox{i}(e_2+e_3),0), ((r \hbox{i}-1) e_1,0), (e_1,e_2+e_3) \rangle ,$

\medskip
\noindent
${\bf h}_4= \langle (e_2+e_3,0), (\hbox{i}(e_2+e_3),0), ((r \hbox{i}-1) e_1,0), (\hbox{i}  e_1,e_2+e_3) \rangle ,$

\medskip
\noindent
${\bf h}_5= \langle (e_2+e_3,0), (\hbox{i}(e_2+e_3),0), ((r \hbox{i}-1) e_1,0), (e_1,e_3) \rangle ,$

\medskip
\noindent
${\bf h}_6= \langle (e_2+e_3,0), (\hbox{i}(e_2+e_3),0), ((r \hbox{i}-1) e_1,0), (\hbox{i} e_1,e_3) \rangle ,$

\medskip
\noindent
${\bf h}_{7}= \langle (e_2+e_3,0), (\hbox{i}(e_2+e_3),0), (\hbox{i} e_1,0), (e_1,e_2+e_3) \rangle $, 

\medskip
\noindent
${\bf h}_{8}= \langle (e_2+e_3,0), (\hbox{i}(e_2+e_3),0), (\hbox{i} e_1,0), (e_1,e_3) \rangle $.

\medskip
\noindent
The exponential image of the subspace ${\bf m}_5$ consists of elements 

\medskip
\noindent
\centerline{$\left ( \left ( \begin{array}{cc} 
a_1+a_2 & b_1+b_2 \hbox{i} \\
b_1-b_2 \hbox{i} & a_1-a_2 \end{array} \right ), \left ( \begin{array}{cc}
c+d & f \\
f & c-d \end{array} \right ) \right ); a_1 \ge 1, c \ge 1, a_2,b_1,b_2,d,f \in \mathbb R $}
   
\medskip
\noindent
with 
$a_1^2-a_2^2-b_1^2-b_2^2=1=c^2-d^2-f^2$. 
Since $p_2(H_i)$, $i=1,2,3,4,7,$ and $p_2(\exp {\bf m}_5)$ satisfy  the 
conditions of 
Lemma 7 the subalgebras $H_i$ for $i=1,2,3,4,7$ cannot  occur as the 
stabilizer of identity for a Bol loop.  
\newline 
The first component of $\exp {\bf h}_i$, $i=5,6,8,$ has the form 
$\left ( \begin{array}{cc} 
\exp \ v & z \\ 
0 & \exp \ -v \end{array} \right )$, where $z \in \mathbb C$, 
 and 
$v=(r \hbox{i}-1)x+ \varepsilon y$, $x,y \in \mathbb R$ with 
$\varepsilon =1$ 
for $j=5$ and  $\varepsilon =i$ in  the case $j=6$, whereas  
$v=\hbox{i} x+y$ for  $j=8$.   

\smallskip
\noindent
The cosets $\left( \left( 
\begin{array}{cc} 
5 & 1 \\
4 & 1 \end{array} \right), 1 \right) H_i$, $i=5,6,8,$ contain the different 
elements $m_1=\left( \left( 
\begin{array}{cc}
\frac{5}{4} & 1 \\
1   & \frac{8}{5} \end{array} \right),1 \right) $ and 
$m_2=\left( \left( \begin{array}{cc}
5 & 4 \\
4   & \frac{17}{5} \end{array} \right),1 \right) $ of $\exp {\bf m}_5$ which is a contradiction to Lemma 7. 
\qed

Author's address: \'A. Figula \\
Mathematisches Institut der Universit\"at Erlangen-N\"urnberg, \\
Bismarkstr. 1 $\frac{1}{2}$,  D-91054 Erlangen, Germany \\
figula@mi.uni-erlangen.de \\
Institute of Mathematics, University of Debrecen, P.O.B. 12, \\
H-4010 Debrecen, Hungary. figula@math.klte.hu

\end{document}